# REPRESENTATION-THEORETIC PROOF OF THE INNER PRODUCT AND SYMMETRY IDENTITIES FOR MACDONALD'S POLYNOMIALS

PAVEL I. ETINGOF[1], ALEXANDER A. KIRILLOV, JR.[2]

## 0. INTRODUCTION

This paper is a continuation of our papers [EK1, EK2]. In [EK2] we showed that for the root system $A_{n-1}$ one can obtain Macdonald's polynomials – a new interesting class of symmetric functions recently defined by I. Macdonald [M1] – as weighted traces of intertwining operators between certain finite-dimensional representations of $U_q\mathfrak{sl}_n$. The main goal of the present paper is to use this construction to give a representation-theoretic proof of Macdonald's inner product and symmetry identities for the root system $A_{n-1}$. Macdonald's inner product identities (see [M2]) have been proved by combinatorial methods by Macdonald (unpublished) for the root system $A_{n-1}$ and by Cherednik in the general case; symmetry identities for the root system $A_{n-1}$ have also been proved by Macdonald ([Macdonald, private communication]).

The paper is organized as follows. In Section 1 we briefly list the basic definitions. In Section 2 we define Macdonald's polynomials $P_\lambda$ and recall the construction of $P_\lambda$ and the inner product between them for the root system $A_{n-1}$ in terms of intertwining operators. By definition, $\langle P_\lambda, P_\mu \rangle = 0$ if $\lambda \neq \mu$, and we show that $\langle P_\lambda, P_\lambda \rangle$ can be expressed as a certain matrix element of product of two intertwining operators. In Section 3 we use the Shapovalov determinant formula to analyze the poles of matrix coefficients of an intertwining operator, and this allows us to express the product of two intertwining operators in terms of a single intertwiner. Applying this to the formula for $\langle P_\lambda, P_\lambda \rangle$ obtained in Section 2, we prove the Macdonald's inner product identity, and the right-hand side is obtained as a product of linear factors in Shapovalov determinant formula. In Section 4 we prove the symmetry identity, which relates the values of $P_\lambda(q^{2(\mu+k\rho)})$ and $P_\mu(q^{2(\lambda+k\rho)})$; the proof is based on the construction of Section 2 and the technique of representing identities in the category of representations of a quantum group by ribbon graphs ([RT1, RT2]). In Section 5 we use the symmetry identities and the fact that Macdonald polynomials are eigenfunctions of certain difference operators (Macdonald operators) to derive recurrence relation for Macdonald polynomials.

### Acknowledgments

We would like to thank our advisor Igor Frenkel for systematic stimulating discussions. We also want to thank Ian Macdonald for fruitful discussions during his visit to Yale in the Fall of 1993; in particular, we learned from him the formulation of the symmetry identities. Finally, this paper was started during our stay at University of Amsterdam as participants of the special period on special functions and representations theory. We would like to express our gratitude to the organizers, especially Tom Koornwinder, for their warm hospitality.

---

[1] Dept. of Mathematics, Harvard Univ., Cambridge, MA 01238; e-mail etingof@math.harvard.edu
[2] Dept. of Mathematics, Yale Univ., New Haven, CT 06520-8283, e-mail:kirillov@math.yale.edu

The work of P.E. was partially supported by the NSF postdoctoral fellowship; the work of A.K. was supported by the Alfred P. Sloan dissertation fellowship.

## 1. Basic definitions

We adopt the following conventions: $\mathfrak{g}$ is a simple Lie algebra over $\mathbb{C}$ of rank $r$, $\mathfrak{h} \subset \mathfrak{g}$ is its Cartan subalgebra, $R \subset \mathfrak{h}^*$ is the corresponding root system, $R^+$ is subset of positive roots, $\alpha_1, \ldots, \alpha_r \in R^+$ is the basis of simple roots, $\theta$ is the highest root. We also introduce the root lattice $Q = \bigoplus \mathbb{Z}\alpha_i$ and the cone of positive roots $Q^+ = \bigoplus \mathbb{Z}_+\alpha_i$.

We fix an invariant bilinear form ( , ) on $\mathfrak{g}$ by the condition that for the associated bilinear form on $\mathfrak{h}^*$ we have $d_i = (\alpha_i, \alpha_i)/2 \in \mathbb{Z}_+$, g.c.d. $(d_i) = 1$; this form allows us to identify $\mathfrak{h}^* \simeq \mathfrak{h} : \lambda \mapsto h_\lambda$. Abusing the language, we will often write, say, $q^\lambda$ instead of $q^{h_\lambda}$.

For every $\alpha \in R$ we define the dual root $\alpha^\vee = \frac{2h_\alpha}{(\alpha,\alpha)} \in \mathfrak{h}$. Let $P = \{\lambda \in \mathfrak{h}^* | \langle \lambda, \alpha^\vee \rangle \in \mathbb{Z}\}$ be the weight lattice, and $P^+ = \{\lambda \in \mathfrak{h}^* | \langle \lambda, \alpha_i^\vee \rangle \in \mathbb{Z}_+\}$ be the cone of dominant weights. Let $\rho = \frac{1}{2}\sum_{\alpha \in R^+} \alpha$; then $\langle \rho, \alpha_i^\vee \rangle = 1$ and thus $\rho \in P^+$. Denote by $W$ the Weyl group for the root system $R$ and by $\mathbb{C}[P]$ the group algebra of the weight lattice, which is spanned by the formal exponentials $e^\lambda, \lambda \in P$. Then $W$ naturally acts on $P$ and on $\mathbb{C}[P]$. Note that in every $W$-orbit in $P$ there is precisely one dominant weight; this implies that the orbitsums $m_\lambda = \sum_{\mu \in W\lambda} e^\mu, \lambda \in P^+$ form a basis of $\mathbb{C}[P]^W$. Finally, we can introduce partial order in $P$: we let $\lambda \leq \mu$ if $\mu - \lambda \in Q^+$.

Let $U_q\mathfrak{g}$ be the quantum group corresponding to $\mathfrak{g}$ (see [D, J] for definitions). We'll use precisely the same form of $U_q\mathfrak{g}$ as we did in [EK2] for $\mathfrak{gl}_n$; in particular, we always consider $q$ as a formal variable and consider $U_q\mathfrak{g}$ and its representations as vector spaces over $\mathbb{C}_q = \mathbb{C}(q^{1/2})$ (we will need the square root of $q$ because of the form of comultiplication we are using). We'll also use the following notions which have been discussed in more detail in [EK2].

We define a polarization of $U_q\mathfrak{g}$ in the usual way: $U_q\mathfrak{g} = U^+ \cdot U^0 \cdot U^-$. For every $\lambda \in \mathfrak{h}^*$ we denote by $M_\lambda$ the Verma module over $U_q\mathfrak{g}$ and by $L_\lambda$ the corresponding irreducible highest-weight module. If $\lambda \in P^+$ then $L_\lambda$ is finite-dimensional. All highest weight modules over $U_q\mathfrak{g}$ have weight decomposition; we will write $V[\alpha]$ for the subspace of homogeneous vectors of weight $\alpha \in \mathfrak{h}^*$ in $V$.

For every finite-dimensional representation $V$ of $U_q\mathfrak{g}$ we define the action of $U_q\mathfrak{g}$ on the space $V^*$ of linear functionals on $V$ by the rule $\langle xv^*, v \rangle = \langle v^*, S(x)v \rangle$ for $v \in V, v^* \in V^*, x \in U_q\mathfrak{g}$, where $S$ is the antipode in $U_q\mathfrak{g}$. This endows $V^*$ with the structure of a $U_q\mathfrak{g}$ representation which we will call the right dual to $V$. In a similar way, the left dual $^*V$ is the representation of $U_q\mathfrak{g}$ in the space of linear functionals on $V$ defined by $\langle xv^*, v \rangle = \langle v^*, S^{-1}(x)v \rangle$. Then the following natural pairings and embeddings are $U_q\mathfrak{g}$-homorphisms:

$$\begin{aligned}
V \otimes {}^*V \to \mathbb{C}, \quad V^* \otimes V \to \mathbb{C} \\
\mathbb{C} \to V \otimes V^*, \quad \mathbb{C} \to {}^*V \otimes V
\end{aligned} \tag{1.1}$$

Note that $V^*$ and ${}^*V$, considered as two structures of a representation of $U_q\mathfrak{g}$ on the same vector space, don't coincide, but they are isomorphic. Namely, $q^{-2\rho}: {}^*V \to V^*$ is an isomorphism. Note also that if $V = L_\lambda$ is an irreducible finite-dimensional representation, so is $V^*$: $L_\lambda^* \simeq L_{\lambda^*}$, where $\lambda^* = -w_0(\lambda)$, $w_0$ being the longest element in the Weyl group.

It is known that if $V, W$ are finite-dimensional then the representations $V \otimes W$ and $W \otimes V$ are isomorphic, but the isomorphism is non-trivial. More precisely (see [D1]), there exists a universal R-matrix $\mathcal{R} \in U_q\mathfrak{g}\hat{\otimes}U_q\mathfrak{g}$ ($\hat{\otimes}$ should be understood as a completed tensor product) such that

$$\check{R}_{V,W} = P \circ \pi_V \otimes \pi_W(\mathcal{R}): V \otimes W \to W \otimes V \tag{1.2}$$

is an isomorphism of representations. Here $P$ is the transposition: $Pv \otimes w = w \otimes v$. Also, it is known that $\mathcal{R}$ has the following form:

$$\begin{aligned}
\mathcal{R} &= q^{-\sum h_i \otimes h_i} \mathcal{R}^*, \quad \mathcal{R}^* \in U^+ \hat{\otimes} U^- \\
(\epsilon \otimes 1)(\mathcal{R}^*) &= (1 \otimes \epsilon)(\mathcal{R}^*) = 1 \otimes 1,
\end{aligned} \tag{1.3}$$

where $\epsilon: U_q\mathfrak{g} \to \mathbb{C}$ is the counit, and $h_i$ is an orthonormal basis in $\mathfrak{h}$.

Similarly to the classical case, one can introduce an involutive algebra automorphism $\omega: U_q\mathfrak{g} \to U_q\mathfrak{g}$ which transposes $U^+$ and $U^-$: $\omega(e_i) = -f_i, \omega(f_i) = -e_i, \omega(h) = -h$. This is a coalgebra antiautomorphism. Thus, for every representation $V$ we can define a new representation $V^\omega$ of $U_q\mathfrak{g}$ in the same space by the formula $\pi_{V^\omega}(x) = \pi_V(\omega x)$. If $v \in V$, we will write $v^\omega$ for the same vector considered as an element of $V^\omega$. If $V$ is finite-dimensional then $V^\omega \simeq {}^*V$ (though the isomorphism is not canonical); in other words, there exists a non-degenerate pairing $(\cdot, \cdot)_V: V \otimes V \to \mathbb{C}$ such that $(xv, v')_V = (v, \omega S(x)v')_V$, which is called the Shapovalov form. This form is symmetric. If $V = L_\lambda$ is irreducible, we will fix this form by the condition that $(v_\lambda, v_\lambda) = 1$. Note that if $v_i, v^i$ are dual bases in $V$ with respect to Shapovalov form then $\mathbf{1}_\lambda = q^{2(\lambda, \rho)}(q^{-2\rho} \otimes 1)\sum v_i \otimes v^i$ is an invariant vector in $V \otimes V^\omega$ such that $\mathbf{1}_\lambda = v_\lambda \otimes v_\lambda^\omega +$ lower terms (by lower terms we always mean terms of lower weight in the first component).

The involution $\omega$ can be extended to intertwiners: if $\Phi: L_\lambda \to L_\nu \otimes U$ is an intertwiner such that $\Phi(v_\lambda) = v_\nu \otimes u_0 +$ lower terms for some $u_0 \in U$ then we can define the intertwiner $\Phi^\omega = P \circ \Phi: L_\lambda^\omega \to U^\omega \otimes L_\nu^\omega$, where $P$ is transposition of $U^\omega$ and $L_\nu^\omega$. Obviously, $\Phi^\omega(v_\lambda^\omega) = u_0^\omega \otimes v_\nu^\omega +$ lower order terms (note that $v_\lambda^\omega$ is a lowest weight vector in $L_\lambda^\omega$).

Finally, we will use the technique of representing homomorphisms in the category of finite-dimensional representations of $U_q\mathfrak{g}$ by ribbon graphs, developed in [RT1, RT2]. For the sake of completeness we briefly recall the basics of this technique in the Appendix.

## 2. Macdonald's polynomials and inner product

Let us briefly recall the definition of Macdonald's polynomials and their construction for root system $A_{n-1}$ in terms of intertwining operators, following the paper [EK2]. Let us fix $k \in \mathbb{Z}_+$.

**Theorem 2.1.** (Macdonald) *There exists a unique family of symmetric trigonometric polynomials $P_\lambda(q, q^k) \in \mathbb{C}(q)[P]^W$ labelled by the dominant weights $\lambda \in P^+$ such that*

1. $P_\lambda(q, q^k) = m_\lambda + \sum_{\mu < \lambda} c_{\lambda\mu} m_\mu$

2. *For fixed $q, k$ the polynomials $P_\lambda(q, q^k)$ are orthogonal with respect to the inner product given by $\langle f, g \rangle_k = \frac{1}{|W|}[f\bar{g}\Delta_{q,q^k}]_0$, where the bar conjugation is defined by $\overline{e^\lambda} = e^{-\lambda}$, $[\ ]_0$ is the constant term of a trigonometric polynomial (i.e., coefficient at $e^0$), and*

$$\Delta_{q,q^k} = \prod_{\alpha \in R} \prod_{i=0}^{k-1}(1 - q^{2i}e^\alpha) \tag{2.1}$$

In our paper [EK2] we showed how these polynomials for the root system $A_{n-1}$ can be obtained from the representation theory of $U_q\mathfrak{sl}_n$. We briefly repeat the main steps here.

From now till the end of this section, we consider only the case $\mathfrak{g} = \mathfrak{sl}_n$. Consider intertwining operators

$$\Phi_\lambda : L_{\lambda + (k-1)\rho} \to L_{\lambda + (k-1)\rho} \otimes U_{k-1}$$

where $U_{k-1} = S^{(k-1)n}\mathbb{C}^n$ is the $q$-deformation of symmetric power of fundamental representation of $U_q\mathfrak{sl}_n$; it can be realized in the space of homogeneous polynomials of degree $n(k-1)$ in $x_1, \ldots, x_n$ (see formula (3.3) below). If $\lambda \in P^+$ then such an intertwiner exists and is unique up to a constant factor. We fix it by the condition $\Phi_\lambda(v_{\lambda^k}) = v_{\lambda^k} \otimes u_0^{k-1} + \ldots$, where $u_0^{k-1} = (x_1 \ldots x_n)^{k-1} \in U_{k-1}[0]$ and $\lambda^k = \lambda + (k-1)\rho$. Define the corresponding "generalized character"

$$\chi_\lambda = \sum_\mu e^\mu \mathrm{Tr}|_{L[\mu]} \Phi_\lambda$$

This is an element of $\mathbb{C}_q[P] \otimes U_{k-1}[0]$. Since $U_{k-1}[0]$ is one-dimensional, it can be identified with $\mathbb{C}$ so that $u_0^{k-1} \mapsto 1$, and thus, $\chi_\lambda$ can be considered as a complex-valued polynomial. Sometimes we will symbolically write $\chi = \mathrm{Tr}(\Phi e^h)$ as an abbreviation of the formula above.

**Theorem 2.2.** ([EK2])

1. 
$$\chi_0 = \prod_{\alpha \in R^+} \prod_{i=1}^{k-1}(e^{\alpha/2} - q^{2i}e^{-\alpha/2})$$

2. $\chi_\lambda$ is divisible by $\chi_0$, and the ratio is a symmetric polynomial.

*3. $\chi_\lambda/\chi_0$ is the Macdonald's polynomial $P_\lambda(q,q^k)$.*

Our main goal is to find an explicit formula for $\langle P_\lambda, P_\lambda\rangle_k$. To do it, note first that it follows from part 1 of the theorem above that $\chi_0\bar\chi_0\Delta_{q,q} = \Delta_{q,q^k}$ and thus

$$(2.2) \qquad \langle P_\lambda, P_\lambda\rangle_k = \frac{1}{|W|}[\chi_\lambda\bar\chi_\lambda\Delta_{q,q}]_0 = \langle \chi_\lambda, \chi_\lambda\rangle_1.$$

Now, it is easy to see from the definition that $\chi_\lambda\bar\chi_\lambda = \text{Tr}(\Psi e^h)$, where $\Psi$ is the following composition

$$(2.3) \quad L_{\lambda^k}\otimes L^\omega_{\lambda^k} \xrightarrow{\Phi_\lambda\otimes\Phi^\omega_\lambda} L_{\lambda^k}\otimes U_{k-1}\otimes U^\omega_{k-1}\otimes L^\omega_{\lambda^k} \xrightarrow{\text{Id}\otimes(\cdot,\cdot)_U\otimes\text{Id}} L_{\lambda^k}\otimes L^\omega_{\lambda^k},$$

where, as before, $\lambda^k = \lambda + (k-1)\rho$. Since the module $L_{\lambda^k}\otimes L^\omega_{\lambda^k}$ is completely reducible, $\text{Tr}(\Psi e^h)$ is just a linear combination of the usual characters of irreducible components in $L_{\lambda^k}\otimes L^\omega_{\lambda^k}$. On the other hand, since the characters for the quantum group are the same as for the classical Lie algebra, and $\Delta_{q,q} = \prod_{\alpha\in R}(1-e^\alpha)$ does not depend on $q$, we know that $\frac{1}{|W|}[(\text{ch }L_\mu)\Delta_{q,q}]_0 = \delta_{\mu,0}$, where ch $L$ is the character of the module $L$ considered as an element of $\mathbb{C}[P]$ (this can be used to prove the orthogonality of $P_\lambda$). Thus, $\langle \chi_\lambda, \chi_\lambda\rangle$ equals the eigenvalue of $\Psi$ on the component in $L_{\lambda^k}\otimes L^\omega_{\lambda^k}$ which is isomorphic to the trivial representation $\mathbb{C}_q$. This gives the following lemma.

**Lemma 2.3.** *The value of the inner product $A_\lambda = \langle P_\lambda, P_\lambda\rangle_k$ can be calculated from the identity*

$$(2.4) \qquad \Psi \mathbf{1}_\lambda = A_\lambda \mathbf{1}_\lambda,$$

*where $\Psi$ is defined by (2.3) and $\mathbf{1}_\lambda$ is an invariant vector in $L_{\lambda^k}\otimes L^\omega_{\lambda^k}$.*

**Theorem 2.4.**

$$(2.5) \qquad \langle P_\lambda, P_\lambda\rangle_k = (\langle v^*_{\lambda^k}, \Phi_\lambda\Phi^\circ_\lambda v_{\lambda^k}\rangle)_{U_{k-1}},$$

*where, as before, $\lambda^k = \lambda + (k-1)\rho$, $(\ )_{U_{k-1}} : U_{k-1}\otimes U^\omega_{k-1} \to \mathbb{C}_q$ is the map given by the Shapovalov form in $U_{k-1}$, and the intertwiner $\Phi^\circ_\lambda : L_{\lambda^k} \to L_{\lambda^k}\otimes U^\omega_{k-1}$ is defined by the condition $\Phi^\circ_\lambda(v_{\lambda^k}) = v_{\lambda^k}\otimes (u_0^{k-1})^\omega +$ lower order terms.*

*Proof.*

The proof is obvious if we use the technique of ribbon graphs. Namely, it follows from Lemma 2.3 that the inner product $\langle P_\lambda, P_\lambda\rangle_k = A_\lambda$ can be defined from the following identity of ribbon graphs:

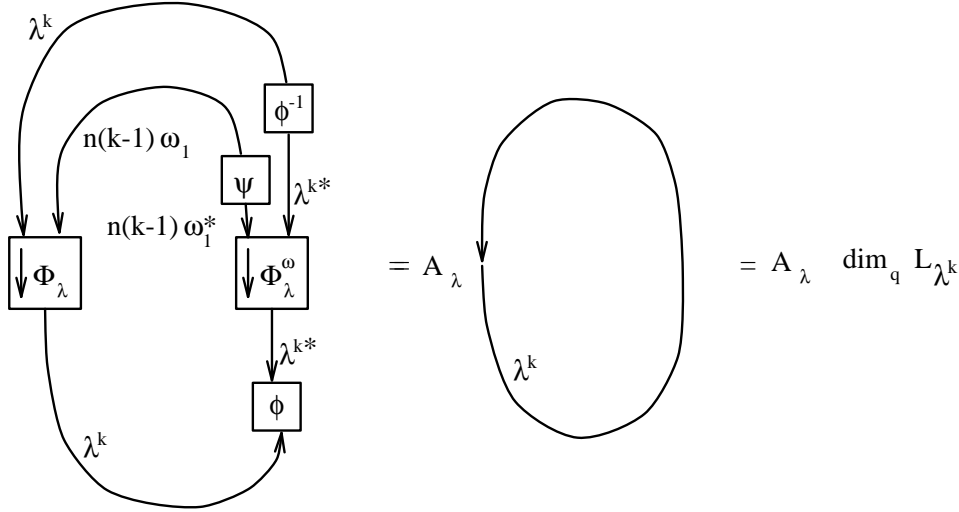

where $\dim_q L = \mathrm{Tr}_L(q^{-2\rho})$, $\phi: L^*_{\lambda^k} \to L_{\lambda^{k*}}$, $\psi: L_{n(k-1)\omega_1^*} \to L^*_{n(k-1)\omega_1} = U^*_{k-1}$ are isomorphisms and $\psi$ is chosen so that $\langle u_0^{k-1}, \psi(u_0^{k-1})^\omega \rangle = 1$. It is easy to check that

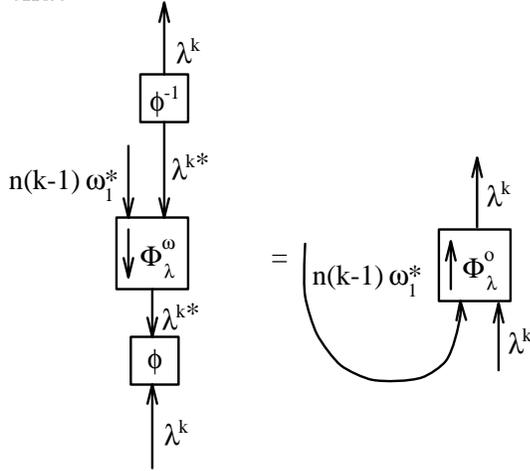

and thus,

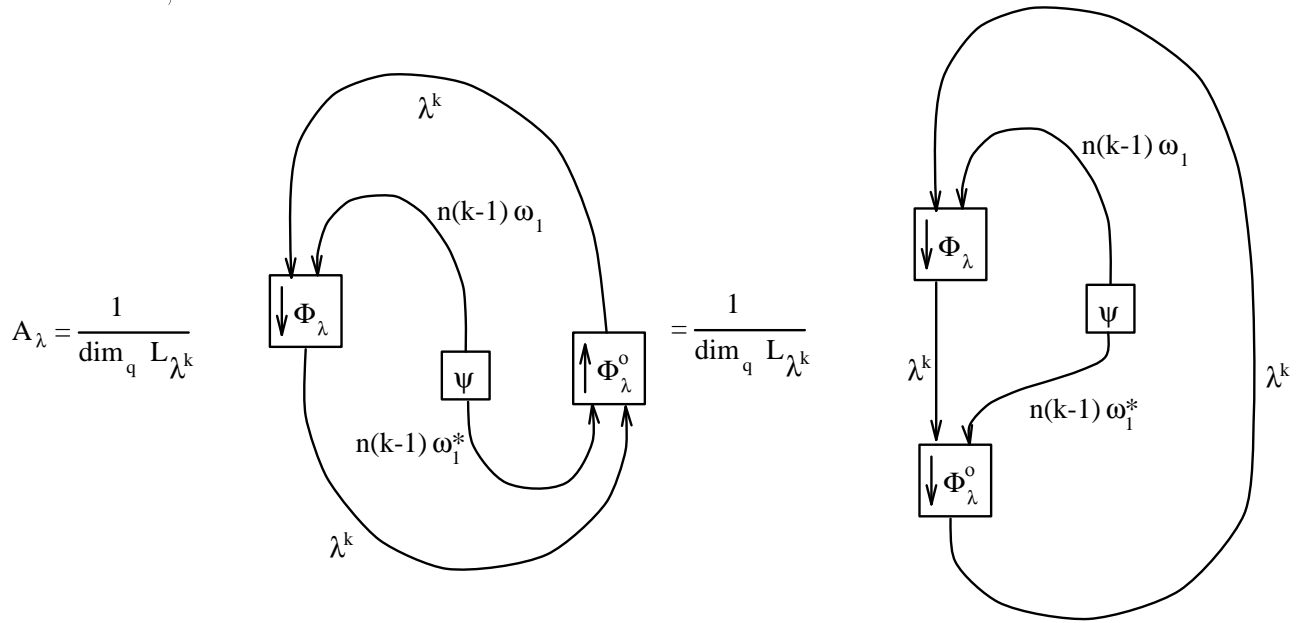

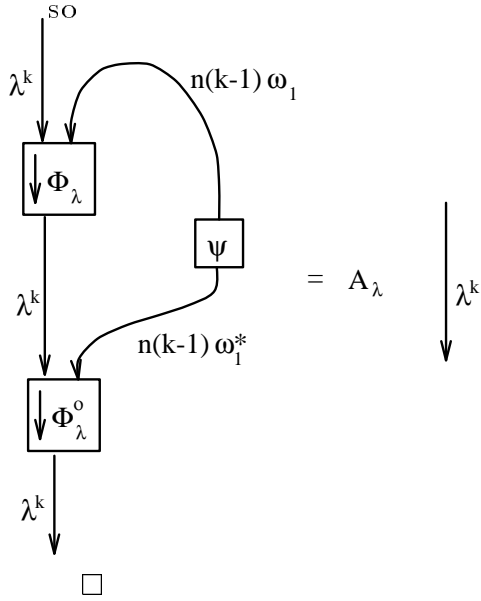

## 3. Algebra of intertwiners and the inner product identity

In this section we consider intertwiners of the form

$$\Phi^\mu_\lambda : M_\lambda \to M_\lambda \otimes L_\mu \tag{3.1}$$

where $M_\lambda$ is the Verma module over $U_q\mathfrak{g}$, $L_\mu$ is the finite-dimensional irreducible module; we assume that $\mu \in P^+ \cap Q$, so $L_\mu[0] \neq 0$. Let $u \in L_\mu[0]$. We consider all the modules over the field of rational functions $\mathbb{C}_q = C(q^{1/2})$; if $\lambda$ is not an integral weight then we also have to add $q^{\langle \lambda, \alpha_i^\vee \rangle/2}$ to this field.

It is known that if $M_\lambda$ is irreducible then there exists a unique intertwiner of the form (3.1) such that $\Phi^\mu_\lambda(v_\lambda) = v_\lambda \otimes u +$ lower order terms. We will denote this intertwiner by $\Phi^{\mu,u}_\lambda$. The same is true if we consider the weight $\lambda$ as indeterminate, i.e. if we consider $t_i = q^{\langle \lambda, \alpha_i^\vee \rangle/2}$ as algebraically independent variables over $\mathbb{C}_q$.

Let us identify $M_\lambda$ with $U^-$ in a standard way. Then we can say that we have a family of actions of $U_q\mathfrak{g}$ in the same space $M \simeq U^-$, and thus we have a family of intertwiners $\Phi^{\mu,u}_\lambda : M \to M \otimes L_\mu$, defined for generic values of $\lambda$.

For $\lambda \in \mathfrak{h}^*$, let us call a trigonometric rational function of $\lambda$ a rational function in $q^{1/2}, q^{\lambda/2}$ (that is, in $q^{1/2}$ and $t_i = q^{\langle \alpha_i^\vee, \lambda \rangle/2}$, $i = 1, \ldots, n$) and call a trigonometric polynomial in $\lambda$ a polynomial in $q^{\pm \lambda/2}$ with coefficients from $\mathbb{C}_q$. Note that the ring of trigonometric polynomials is a unique factorization ring, and invertible elements in this ring are of the form $c(q)q^{\langle \lambda, \alpha \rangle}$, $\alpha \in \frac{1}{2}Q^\vee$.

**Lemma 3.1.** *For fixed $\mu \in P^+, u \in L_\mu[0], u \neq 0$ we have*

1. *For generic $\lambda$, there exists a unique intertwining operator $\Phi^{\mu,u}_\lambda : M_\lambda \to M_\lambda \otimes L_\mu$ such that $\Phi^{\mu,u}_\lambda(v_\lambda) = v_\lambda \otimes u + l.o.t.$. Its matrix elements are trigonometric rational functions of $\lambda$.*

2. *The least common denominator of the matrix elements of $\Phi^{\mu,u}_\lambda$ is given by*

$$(3.2) \qquad d_\mu(\lambda) = \prod_{\alpha \in R^+} \prod_{i=1}^{n_\mu^\alpha} \left(1 - q^{(2(\alpha, \lambda + \rho) - i(\alpha, \alpha))}\right),$$

where $n_\mu^\alpha = \max\{i \in \mathbb{Z}_+ | \mu - i\alpha \in Q^+\}$.

The second part should be understood in the following sense: if we define $\tilde\Phi_\lambda^{\mu,u} = d_\mu(\lambda)\Phi_\lambda^{\mu,u}$ then $\tilde\Phi_\lambda^{\mu,u}$ is an intertwiner $M_\lambda \to M_\lambda \otimes L_\mu$ whose matrix coefficients are trigonometric polynomials in $\lambda$ (thus, it is defined for all $\lambda \in \mathfrak{h}^*$) and $\tilde\Phi_\lambda^{\mu,u} = p(\lambda)\Phi'$ for some intertwiner $\Phi'$ with polynomial coefficients and trigonometric polynomial $p(\lambda)$ is possible only if $p(\lambda)$ is invertible.

*Proof.* Proof of Lemma 3.1 is quite similar to the classical case (cf. [ES]); we briefly sketch it here. First, existence and uniqueness of $\Phi_\lambda^{\mu,u}$ for generic $\lambda$ follows from general arguments: it is known that the space of intertwiners $M_\lambda \to M_\lambda^c \otimes L_\mu$, where $M_\lambda^c$ is the contragredient Verma module, is isomorphic to $L_\mu[0]$; on the other hand, if $M_\lambda$ is irreducible then $M_\lambda \simeq M_\lambda^c$. Since the condition of being an intertwiner can be written as a system of linear equations on components of $\Phi$ with polynomial coefficents, solution can be written as a trigonometric rational function. This proves the first part of the theorem.

Next, one can write explicit formula for $\Phi$: $\Phi(v_\lambda) = \sum_{k,l}(F^{-1})_{kl}g_k v_\lambda \otimes (\omega g_l)u$, where $g_k$ is a basis in $U^-$, and $F^{-1}$ is the inverse matrix to the Shapovalov form (cf. [ES]). Thus, it may have poles only at the points where the determinant of Shapovalov form vanishes. The formula for the determinant of the Shapovalov form in the quantum case can be found, for example, in [CK], and the factors occuring there are precisely the factors in formula (3.2) (up to invertible factors). Finally, we have to show that in fact the poles actually do occur on the hyperplanes which enter formula (3.2) and all the poles are simple; this can be done precisely in the same way as for $q = 1$, i.e. in the classical case (see [ES]).

*Remark.* It is seen from this proof that $\tilde\Phi_\lambda^{\mu,u}$ is actually a trigonometric polynomial in $\lambda$ with operator coefficients, i.e. the degrees of its matrix coefficients, as trigonometric polynomials, are uniformly bounded (under a suitable definition of degree).

We will also need one more technical lemma.

**Lemma 3.2.** *Let us write $\tilde\Phi_\lambda^{\mu,u}$ in the following form:*

$$\tilde\Phi_\lambda^{\mu,u} v_\lambda = d_\mu(\lambda) v_\lambda \otimes u + \ldots + a(\lambda) v_\lambda \otimes u_\mu,$$

*where $u_\mu$ is the highest-weight vector in $L_\mu$, and $a(\lambda) \in \mathbb{C}_q[q^{\pm\lambda/2}] \otimes U^-[-\mu]$ is a trigonometric polynomial of $\lambda$ with values in the universal enveloping algebra. Then the greatest common divisor of the components of $a(\lambda)$ is 1.*

*Proof.* It is easy to see, using the irreducibility of $L_\mu$, that if $a(\lambda) = 0$ then $\tilde\Phi_\lambda^{\mu,u} = 0$. On the other hand, we have shown before that the coefficients of $\tilde\Phi_\lambda^{\mu,u}$ have no nontrivial common divisors, and thus $\tilde\Phi_\lambda^{\mu,u}$ could only vanish on a subvariety of codimension more than one. Thus, the same must be true for $a(\lambda)$.

Now we want to define a structure of algebra on these intertwiners. Let $\Phi_1 : M_\lambda \to M_\lambda \otimes L_{\mu_1}, \Phi_2 : M_\lambda \to M_\lambda \otimes L_{\mu_2}$ be non-zero intertwiners. Let us define their product $\Phi_1 * \Phi_2 : M_\lambda \to M_\lambda \otimes L_{\mu_1+\mu_2}$ as the composition

$$M_\lambda \xrightarrow{\Phi_2} M_\lambda \otimes L_{\mu_2} \xrightarrow{\Phi_1 \otimes 1} M_\lambda \otimes L_{\mu_1} \otimes L_{\mu_2} \xrightarrow{1 \otimes \pi} M_\lambda \otimes L_{\mu_1+\mu_2}$$

where $\pi$ is a fixed projection $\pi : L_{\mu_1} \otimes L_{\mu_2} \to L_{\mu_1+\mu_2}$.

Let us now consider a very special case of the above situation. *From this moment till the end of this section we only work with $\mathfrak{g} = \mathfrak{sl}_n$.* Take $\mu = kn\omega_1$ for some $k \in \mathbb{Z}_+$, that is, $L_\mu = S^{kn}\mathbb{C}^n$, where $\mathbb{C}^n$ is the fundamental representation of $U_q\mathfrak{sl}_n$. One can easily write the action of $U_q\mathfrak{sl}_n$ in this space in terms of difference derivatives (see [EK2]):

(3.3)
$$h_i \mapsto x_i \frac{\partial}{\partial x_i} - (k-1), \ e_i \mapsto x_i D_{i+1}, \ f_i \mapsto x_{i+1} D_i$$
$$(D_i f)(x_1, \ldots, x_n) = \frac{f(x_1, \ldots, qx_i, \ldots, x_n) - f(x_1, \ldots, q^{-1}x_i, \ldots, x_n)}{(q - q^{-1})x_i}$$

In this case all the weight subspaces of $L_\mu$ are one-dimensional; in particular, we can choose $u_0^k = (x_1 \ldots x_n)^k \in L_\mu[0]$, where $x_1, \ldots, x_n$ is the basis in $\mathbb{C}^n$; then $L_\mu[0] = \mathbb{C}u_0^k$. For brevity, we will write $U_k$ for $L_{kn\omega_1}$, $\Phi_\lambda^k$ for $\Phi_\lambda^{\mu=kn\omega_1, u_0^k}$, etc. Let us fix the projection $\pi : U_k \otimes U_l \to U_{k+l}$ by $\pi(u_0^k \otimes u_0^l) = u_0^{k+l}$. In this case, $n_\mu^\alpha = k$ for all $\alpha \in R^+$ and

(3.4)
$$d_k(\lambda) = \prod_{\alpha \in R^+} \prod_{i=1}^{k} (1 - q^{2(\alpha, \lambda+\rho)-2i})$$

**Theorem 3.3.**

(3.5)
$$\tilde{\Phi}_\lambda^k * \tilde{\Phi}_\lambda^l = \tilde{\Phi}_\lambda^{k+l}$$

*Proof.* Let us denote the left-hand side of (3.5) by $\Psi$. Then $\Psi$ is an intertwiner $M_\lambda \to M_\lambda \otimes U_{k+l}$, whose matrix coefficients are trigonometric polynomials in $\lambda$. In particular, we can write $\Psi(v_\lambda) = f(\lambda)v_\lambda \otimes u_0^{k+l} + $ l.o.t. On the other hand $\tilde{\Phi}_\lambda^{k+l}(v_\lambda) = d_{k+l}(\lambda)v_\lambda \otimes u_0^{k+l} + $ l.o.t. Since the intertwining operator is unique for generic $\lambda$, this implies $\Psi(\lambda) = \frac{f(\lambda)}{d_{k+l}(\lambda)}\tilde{\Phi}_\lambda^{k+l}$. Since the greatest common divisor of the matrix elements of $\tilde{\Phi}^{k+l}$ is 1, this implies that $f(\lambda)$ is divisible by $d_{k+l}(\lambda)$.

Let us now consider the lowest term of $\Psi$. If we write the lowest term of $\tilde{\Phi}_\lambda^k$ as $a_k(\lambda)v_\lambda \otimes u_k$ (cf. Lemma 3.3) and lowest term of $\Phi_\lambda^l$ as $a_l(\lambda)v_\lambda \otimes u_l$ then the lowest term of $\Psi$ will be $a_l(\lambda)a_k(\lambda)v_\lambda \otimes u_{k+l}$ (up to some power of $q$). Since we know that components of $a_k$ have no common divisors, and the same is true for $a_l$, it follows that the greatest common divisor of components of $a_k(\lambda)a_l(\lambda)$ is 1. Indeed, suppose that $p(\lambda)$ is a common divisor of components of $a_k(\lambda)a_l(\lambda)$. Passing if necessary to a certain algebraic extension of $\mathbb{C}_q$ we get that $a_k(\lambda)a_l(\lambda)$ vanishes on a certain subvariety of codimension 1. On the other hand, this contradicts to the fact that both $a_k, a_l$ could only vanish on subvarieties of codimension more than one, since $U_q\mathfrak{g}$ has no zero divisors. Thus, the greatest common divisor of coefficients of $\Psi(\lambda)$ is one, which implies that $d_{k+l}(\lambda)$ is divisible by $f(\lambda)$.

This proves that $\Psi(\lambda) = c(q)q^{(\lambda,\alpha)}\tilde{\Phi}_\lambda^{k+l}$ for some $\alpha \in \frac{1}{2}Q$ and rational function $c(q)$, independent of $\lambda$. To calculate $\alpha, c$, let us introduce the notion of degree: $\deg(q^{(\lambda,\alpha)}) = -(\alpha, \rho)$ and consider the terms of maximal degree on both sides of (3.5). Heuristically, this corresponds to the asymptotics as $\lambda \to +\infty\rho$.

**Lemma 3.4.** *If one denotes the maximal degree term of $\tilde{\Phi}^k$ by $\bar{\Phi}^k$ then*

$$\bar{\Phi}_k(v_\lambda) = v_\lambda \otimes u_0^k.$$

To prove the lemma, note first that the maximal degree term of $\tilde{\Phi}^k$ is equal to the maximal degree term of $\Phi^k$. Now, write $\Phi(v_\lambda) = \sum_{k,l}(F^{-1})_{kl}g_k v_\lambda \otimes (\omega g_l)u$. It is known that if we choose a basis $g_k$ in such a way that $g_0 = 1, g_k$ has strictly negative weight for $k > 0$, then in the limit $\lambda \to \infty\rho$ one has $(F^{-1})_{kl} = \delta_{k,0}\delta_{l,0}$ (this, for example, follows from [L, Proposition 19.3.7], which gives much more detailed information about the asymptotic behavior of $F$; it states that under a suitable normalization the Shapovalov form in the Verma module $M_\lambda$ formally converges to the Drinfeld's form on $U_q(\mathfrak{n}^-)$ as $\lambda \to +\infty\rho$).

Using this lemma and the fact that in the identification $M_\lambda \simeq M \simeq U^-$ the action of $U^-$ does not depend on $\lambda$, one can show that if we denote the highest degree term of $\Psi$ by $\bar{\Psi}$ then

$$\bar{\Psi}(v_\lambda) = v_\lambda \otimes u_{k+l}.$$

Comparing it with the the expression for $\bar{\Phi}_{k+l}(v_\lambda)$, we get the statement of the theorem. □

**Corollary 3.5.**

$$(3.6) \qquad \Phi_\lambda^{k+l} = \frac{d_k(\lambda)d_l(\lambda)}{d_{k+l}(\lambda)}\Phi_\lambda^k * \Phi_\lambda^l$$

So far, we have proved Theorem 3.3 only for the case when $k, l \in Z_+$. However, it can be generalized. Let us consider the space $\tilde{U}_k = \{(x_1 \ldots x_n)^k p(x), p(x) \in \mathbb{C}_q[x_1^{\pm 1}, \ldots x_n^{\pm 1}], p(x)$ is a homogeneous polynomial of degree 0$\}$ where $k$ is an arbitrary complex number. Formula (3.3) defines an action of $U_q\mathfrak{sl}_n$ in $\tilde{U}_k$. Also, define $u_0^k = (x_1 \ldots x_n)^k \in \tilde{U}_k$.

**Lemma 3.6.**

1. *The set of weights of $\tilde{U}_k$ coincides with the weight lattice $Q$, and each weight subspace is one-dimensional. In particular, $U_k[0] \simeq \mathbb{C}u_0^k$.*

2. *For generic $k$, the mapping*

$$(3.7) \qquad x^\lambda \mapsto \frac{\tilde{\Gamma}_q(\lambda_1 + 1)\ldots\tilde{\Gamma}_q(\lambda_n + 1)}{(\tilde{\Gamma}_q(k+1))^n}x^{-1-\lambda},$$

defines an isomorphism $\tilde{U}_k^\omega \simeq \tilde{U}_{-1-k}$. The normalization is chosen so that $u_0^k \mapsto u_0^{-1-k}$. Here $\tilde{\Gamma}_q(\lambda)$ is renormalized q-gamma function:

$$\tilde{\Gamma}_q(x) = q^{-(x-1)(x-2)/2} \frac{1}{(1-q^2)^{x-1}} \prod_{n=0}^{\infty} \frac{1-q^{2(n+1)}}{1-q^{2(n+x)}},$$

so $\tilde{\Gamma}(x+1) = [x]\tilde{\Gamma}(x)$, where

$$[x] = \frac{q^x - q^{-x}}{q - q^{-1}}.$$

3. If $k \in \mathbb{Z}_+$ then $\tilde{U}_k$ contains a finite-dimensional submodule, isomorphic to the module $U_k$ defined above: $U_k = \tilde{U}_k \cap \mathbb{C}(q)[x_1, \ldots, x_n]$. Also, in this case $\tilde{U}_{-1-k}$ has a finite-dimensional quotient $U^k = \tilde{U}_{-1-k}/(x^\lambda$ such that at least one $\lambda_i \in \mathbb{Z}_+)$. In particular, $U_{-1}$ can be projected onto $U^0 \simeq \mathbb{C}$. Moreover, formula (3.7) above defines an isomorphism $U_k^\omega \simeq U^k$ for $k \in \mathbb{Z}_+$.

*Proof* of this lemma is straightforward.

Now, let us assume that $\lambda$ is generic and consider an intertwiner $\Phi_\lambda^k : M_\lambda \to M_\lambda \hat{\otimes} \tilde{U}_k$ such that $\Phi_\lambda^k(v_\lambda) = v_\lambda \otimes u_0^k + \ldots$, and $\hat{\otimes}$ is a tensor product completed with respect to $\rho$-grading in $M_\lambda$. Note that if $k \in \mathbb{Z}_+$ then image of $\Phi(v_\lambda)$ lies in the submodule $M_\lambda \otimes U_k$ (which follows, for example, from the explicit formula for $\Phi$, see [ES]), so this is consistent with our previuos notations. Also, for $k \in \mathbb{C}$ we define

(3.8) $$d_k(\lambda) = \prod_{\alpha \in R^+} \prod_{i=0}^{\infty} \frac{1 - q^{2(\alpha, \lambda+\rho)+2i} q^{-2k}}{1 - q^{2(\alpha, \lambda+\rho)+2i}},$$

which we can consider as a formal power series in $q$ with coefficients that are rational functions in $q^\lambda, q^k$ (which we consider as independent variables). Note that if $k \in \mathbb{Z}_+$, this coincides with previously given definition.

**Theorem 3.7.** *For any $k \in \mathbb{Z}_+, l \in \mathbb{C}$ we have*

(3.9) $$\Phi_\lambda^k * \Phi_\lambda^l = \frac{d_{k+l}(\lambda)}{d_k(\lambda) d_l(\lambda)} \Phi_\lambda^{k+l},$$

*where $d_k(\lambda)$ is given by formula* (3.8).

*Proof.* Let us fix $k$. Then the matrix elements of the operators on both sides of (3.9) are rational functions in $q, q^l, q^\lambda$, which follows from the fact that $\frac{d_{k+l}(\lambda)}{d_k(\lambda)d_l(\lambda)}$ is a rational function in $q, q^l, q^\lambda$. Now the statement of the lemma follows from the fact that this is true for $l \in \mathbb{Z}_+$ and the following trivial statement:

If $F(q,t) \in \mathbb{C}(q,t)$ is such that $F(q, q^l) = 0$ for all $l \in \mathbb{Z}_+$ then $F = 0$.

□

Let us apply this to case when $l = -1 - k$. In this case explicit calculation gives the following answer:

**Corollary 3.8.** *For $k \in \mathbb{Z}_+$,*

$$\Phi_\lambda^k * \Phi_\lambda^{-1-k} = \Phi_\lambda^{-1} \prod_{\alpha \in R^+} \prod_{i=1}^{k} \frac{1 - q^{2(\alpha,\lambda+\rho)+2i}}{1 - q^{2(\alpha,\lambda+\rho)-2i}}$$

Now we are in the position to prove Macdonald's inner product identities. Namely, in Section 2 we have proved that

$$\langle P_\lambda, P_\lambda \rangle_k = (\langle v_{\lambda^k}^*, \Phi_1 \Phi_2 v_{\lambda^k} \rangle)_{U_{k-1}},$$

where $\lambda^k = \lambda + (k-1)\rho$, the intertwiner $\Phi_1 : L_{\lambda^k} \to L_{\lambda^k} \otimes U_{k-1}$ is such that $\Phi_1(v_{\lambda^k}) = v_{\lambda^k} \otimes u_0^{k-1} + \ldots$, $\Phi_2 : L_{\lambda^k} \to L_{\lambda^k} \otimes U_{k-1}^\omega$ is such that $\Phi_2(v_{\lambda^k}) = v_{\lambda^k} \otimes (u_0^{k-1})^\omega + \ldots$, and the Shapovalov form $(\,,\,)_{U_{k-1}} : U_{k-1} \otimes U_{k-1}^\omega \to \mathbb{C}$ is normalized so that $(u_0^{k-1}, (u_0^{k-1})^\omega) = 1$ (there is no contradiction with previous conventions, since we have the freedom of fixing the highest weight vector in $U_{k-1}$). Comparing this with the results and notations of this Section, we see that

$$\langle P_\lambda, P_\lambda \rangle_k = \operatorname{Res}(\langle v_{\lambda^k}^*, (\Phi_{\lambda^k}^{k-1} * \Phi_{\lambda^k}^{-k}) v_{\lambda^k} \rangle),$$

and Res is the projection $U_{-1} \to \mathbb{C}$ such that $u_0^{-1} \mapsto 1$. Using Corollary 3.8, we immediately obtain

**Theorem 3.9.** (Macdonald)

$$\langle P_\lambda, P_\lambda \rangle_k = \prod_{\alpha \in R^+} \prod_{i=1}^{k-1} \frac{1 - q^{2(\alpha,\lambda+k\rho)+2i}}{1 - q^{2(\alpha,\lambda+k\rho)-2i}}.$$

This is precisely the Macdonald's inner product identity for the root system $A_{n-1}$.

*Remark.* It is easy to generalize the inner product identity to the case when $k$ is a generic complex number (here we assume that $q$ is specialized to a fixed real number between 0 and 1). The Macdonald's polynomials in this situation are defined in the same way as in Section 2 (Theorem 2.1), but in this case

$$(3.10) \qquad \Delta_{q,q^k} = \prod_{\alpha \in R} \prod_{i=0}^{\infty} \frac{1 - q^{2i} e^\alpha}{1 - q^{2k} q^{2i} e^\alpha},$$

which coincides with (2.1) for positive integer values of $k$. It is easy to show that (3.10) defines a Laurent series in $q$ whose coefficients are analytic functions of $t = q^k$ for $|t| \leq q$. Also, it is easy to see that the coeffecents of Macdonald's polynomials are rational functions of $t = q^k$ which are smooth at $t = 0$. Therefore, they are defined for a generic value of $t$, and hence we have

$$(3.11) \qquad \langle P_\lambda, P_\lambda \rangle_k = \prod_{\alpha \in R^+} \prod_{i=1}^{\infty} \frac{(1 - q^{2(\alpha,\lambda+k\rho)+2i})(1 - q^{2(\alpha,\lambda+k\rho)+2(i-1)})}{(1 - q^{2(\alpha,\lambda+k\rho)+2(i+k-1)})(1 - q^{2(\alpha,\lambda+k\rho)+2(i-k)})}.$$

Indeed, both sides of (3.11) are holomorphic in $t = q^k$ for $|t| < q$, and coincide for $t = q^k$, $k \in \mathbb{Z}^+$ (since for this case (3.11) reduces to Theorem 3.9). Therefore, they must coincide identically.

# 4. Symmetry identity.

In this section we only consider the case $\mathfrak{g} = \mathfrak{sl}_n$.

The main goal of this section is to prove Theorem 4.3, which establishes certain symmetry between the values $P_\lambda(q^{2(\mu+k\rho)})$ and $P_\mu(q^{2(\lambda+k\rho)})$ (notations will be explained later). The proof of this theorem is based on the technique of ribbon graphs.

As in Section 3, let $\Phi_\lambda : L_{\lambda^k} \to L_{\lambda^k} \otimes U_{k-1}$ be such that $\Phi(v_{\lambda^k}) = v_{\lambda^k} \otimes u_0^{k-1} + \ldots$, $\lambda^k = \lambda + (k-1)\rho$.

**Lemma 4.1.**

(4.1)

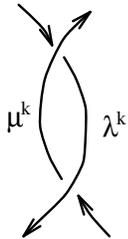

$$= \chi_\mu(q^{2(\lambda+k\rho)})\Phi_\lambda$$

where $\chi_\mu \in \mathbb{C}[P]$ is the weighted trace of $\Phi_\mu$ (see Section 2), and $\chi_\mu(q^\lambda)$ stands for polynomial in $q, q^{-1}$ which is obtained by replacing each formal exponent $e^\alpha$ in the expression for $\chi_\mu$ by $q^{(\alpha,\lambda)}$.

*Proof.* Let us consider the operator $F : L_{\lambda^k} \to L_{\lambda^k} \otimes U_{k-1}$ corresponding to the ribbon graph on the left hand side of (4.1). It is some $U_q\mathfrak{g}$-homomorphism. Since we know that such a homomorphism is unique up to a constant, it follows that $F = a\Phi_\lambda$ for some constant $a$. To find $a$, let us find the image of the highest-weight vector. First, consider the following part of this picture:

Then it follows from the explicit form of R-matrix (1.3) that the corresponding operator $A : V_{\lambda^k} \otimes (V_{\mu^k})^* \to V_{\lambda^k} \otimes (V_{\mu^k})^*$ has the following property: if $x \in V_{\mu^k}^*[\alpha]$ then $A(v_{\lambda^k} \otimes x) = q^{-2(\lambda^k,\alpha)} v_{\lambda^k} \otimes x + \ldots$. Thus, if $x_i$ is basis in $L_{\mu^k}$, $x^i$ – dual basis in $L_{\mu^k}^*$, $x_i$ has weight $\alpha_i$ then explicit calculation shows that

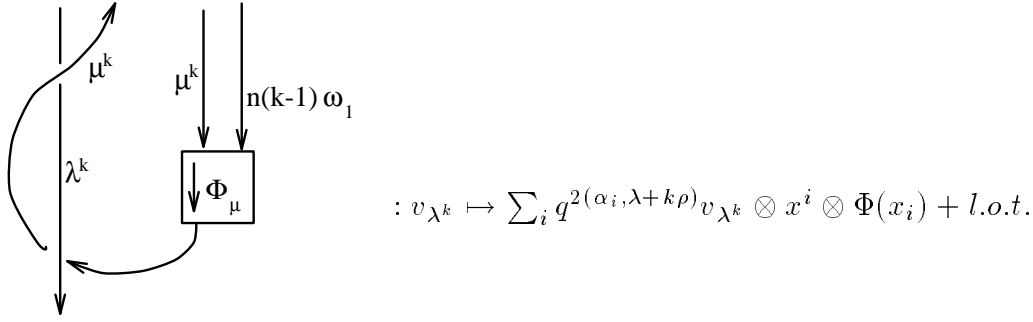

$$: v_{\lambda^k} \mapsto \sum_i q^{2(\alpha_i, \lambda + k\rho)} v_{\lambda^k} \otimes x^i \otimes \Phi(x_i) + l.o.t.$$

and thus, $F(v_{\lambda^k}) = \chi_\mu(q^{2(\lambda+k\rho)}) v_{\lambda^k} \otimes u_0^{k-1} + \ldots$, which completes the proof. □

**Corollary 4.2.** *Let* $\Phi_\lambda : M_\lambda \to M_\lambda \otimes U_{k-1}, \Phi_\lambda^\circ : M_\lambda \to M_\lambda \otimes U_{k-1}^\omega$ *be as in Theorem 2.4. Then*

(4.2)

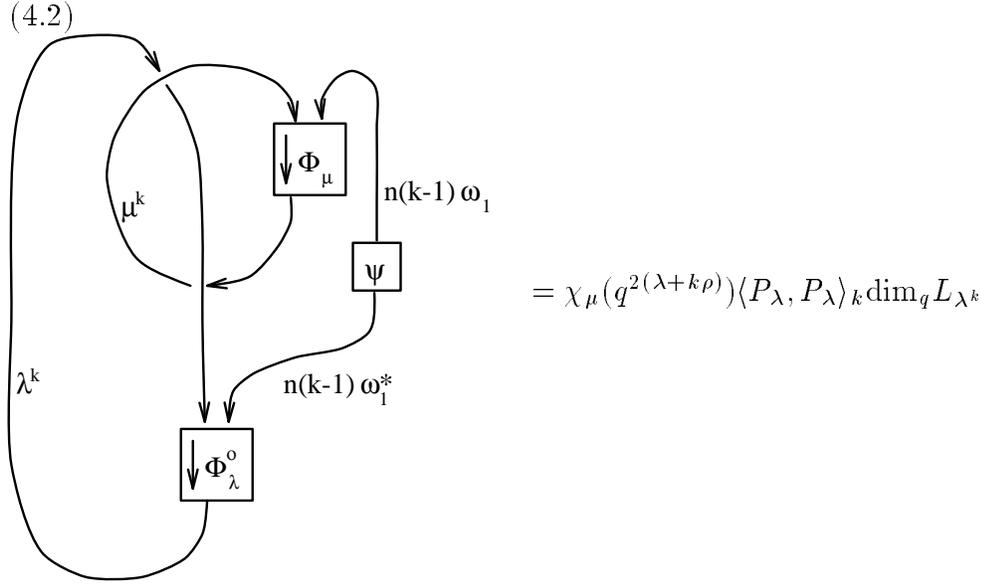

$$= \chi_\mu(q^{2(\lambda+k\rho)}) \langle P_\lambda, P_\lambda \rangle_k \dim_q L_{\lambda^k}$$

*Proof.* This follows from the previous lemma and the arguments used in the proof of Theorem 2.4. □

In a similar way, repeating with necessary changes all the steps of Sections 2 and 3 one can prove

**Corollary 4.2'.** *Formula (4.2) remains valid if we replace in the graph on the right hand side* $\Phi_\mu$ *by* $\Phi_\mu^\circ$, $\Phi_\lambda^\circ$ *by* $\Phi_\lambda$ *and interchange* $\omega_1$ *and* $\omega_1^*$.

**Theorem 4.3.**
(4.3)
$$\frac{P_\mu(q^{2(\lambda+k\rho)})}{P_\lambda(q^{2(\mu+k\rho)})} = q^{2k(\rho, \lambda-\mu)} \prod_{\alpha \in R^+} \prod_{i=0}^{k-1} \frac{1 - q^{2(\alpha, \mu+k\rho)+2i}}{1 - q^{2(\alpha, \lambda+k\rho)+2i}} = \prod_{\alpha \in R^+} \prod_{i=0}^{k-1} \frac{[(\alpha, \mu+k\rho)+i]}{[(\alpha, \lambda+k\rho)+i]}$$

*Proof.* The proof is based on the following identity of the ribbon graphs:

(4.4)

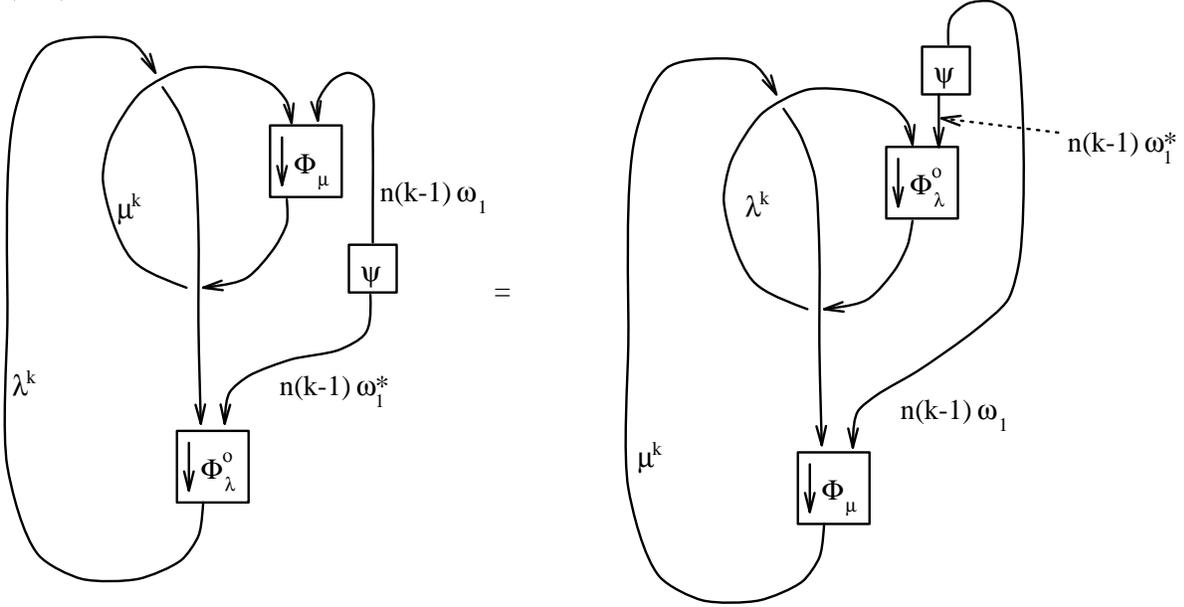

Due to Corollaries 4.2 and 4.2′, this implies:

(4.5) $\qquad \chi_\mu(q^{2(\lambda+k\rho)})\langle P_\lambda, P_\lambda\rangle_k \dim_q L_{\mu^k} = \chi_\lambda(q^{2(\mu+k\rho)})\langle P_\mu, P_\mu\rangle_k \dim_q L_{\lambda^k}$

Substituting in this formula explicit expression for $\langle P_\lambda, P_\lambda\rangle_k$ (Theorem 3.9) and using the expression for $\dim_q L_\lambda$:

$$\dim_q L_\lambda = q^{-2(\lambda,\rho)} \prod_{\alpha\in R^+} \frac{1-q^{2(\alpha,\lambda+\rho)}}{1-q^{2(\alpha,\rho)}} = \prod_{\alpha\in R^+} \frac{[(\alpha,\lambda+\rho)]}{[(\alpha,\rho)]},$$

which can be easily deduced from the Weyl character formula, we get the statement of the theorem. □

**Corollary 4.4.** (Macdonald's special value identity, [M1,M2])
(4.6)
$$P_\lambda(q^{2k\rho}) = q^{-2k(\rho,\lambda)} \prod_{\alpha\in R^+} \prod_{i=0}^{k-1} \frac{1-q^{2(\alpha,\lambda+k\rho)+2i}}{1-q^{2(\alpha,k\rho)+2i}} = \prod_{\alpha\in R^+} \prod_{i=0}^{k-1} \frac{[(\alpha,\lambda+k\rho)+i]}{[(\alpha,k\rho)+i]},$$

*Proof.* Let $\mu=0$. Then $P_\mu = 1$, and formula (4.3) reduces to (4.6).□

*Remark.* Similarly to the arguments at the end of Section 3, we can show that for generic $k$

$$\frac{P_\mu(q^{2(\lambda+k\rho)})}{P_\lambda(q^{2(\mu+k\rho)})} = q^{2k(\rho,\lambda-\mu)} \prod_{\alpha\in R^+} \prod_{i=0}^{\infty} \frac{1-q^{2(\alpha,\mu+k\rho)+2i}}{1-q^{2(\alpha,\lambda+k\rho)+2i}} \frac{1-q^{2(\alpha,\lambda+k\rho)+2(k+i)}}{1-q^{2(\alpha,\mu+k\rho)+2(k+i)}}.$$

Actually, in this case we don't even need an analytic argument: it is enough to observe that both sides of this equality are rational in $q^k$ and coincide for positive integer values of $k$.

## 5. Recursive relations

In this section we explain how recurrence relations for Macdonald's polynomials can be deduced from the symmetry identity.

For $r \in \{1, ..., n-1\}$, let $\Lambda_r$ denote the set of weights of the representation $\Lambda^r \mathbb{C}^n$ of $\mathfrak{sl}_n$. It can be naturally identified with the set of subsets of size $r$ in $\{1, ..., n\}$.

Recall [M1] that Macdonald's polynomials satisfy the difference equations

$$(5.1) \qquad (M_r P_\lambda) = c_\lambda^r P_\lambda,$$

where $c_\lambda^r = \sum_{\nu \in \Lambda_r} q^{2(\lambda + k\rho, \nu)} = X_r(q^{2(\lambda + k\rho)})$, $X_r$ being the character of $\Lambda^r \mathbb{C}^n$, and $M_r$ are the operators

$$(5.2) \qquad M_r = q^{kr(r-n)} \sum_{\nu \in \Lambda_r} \left( \prod_{\alpha \in R: (\alpha, \nu) = -1} \frac{q^{2k} - e^\alpha}{1 - e^\alpha} \right) T_\nu,$$

where $T_\nu e^\lambda = q^{2(\nu, \lambda)} e^\lambda$. (These operators were introduced independently by I. Macdonald and S. Ruijsenaars.)

Let us specialize (5.1) at $q^{2(\mu + k\rho)}$, $\mu \in P^+$. Notice that under this specialization in the right hand side of (5.2) the terms corresponding to such values of $\nu$ that $\mu + \nu \notin P^+$ drop out. Indeed, if $\mu + \nu \notin P^+$ then we are forced to have $(\mu, \alpha_j) = 0$, $(\nu, \alpha_j) = -1$ for some simple root $\alpha_j$. The coefficient to $T_\nu$ in (5.2) is $\prod_{\alpha \in R: (\alpha, \nu) = -1} \frac{q^{2k} - q^{2(\mu + k\rho, \alpha)}}{1 - q^{2(\mu + k\rho, \alpha)}}$. Clearly, it contains the factor $q^{2k} - q^{2(\mu + k\rho, \alpha_j)}$ in the numerator. But this factor is zero since $(\mu, \alpha_j) = 0$. So the whole coefficient is zero and the term drops out. Thus we get

$$\sum_{\nu \in \Lambda_r : \mu + \nu \in P^+} \left( \prod_{\alpha \in R: (\alpha, \nu) = -1} \frac{[(\mu + k\rho, \alpha) - k]}{[(\mu + k\rho, \alpha)]} \right) P_\lambda(q^{2(\mu + \nu + k\rho)})$$
$$(5.3) \qquad = X_r(q^{2(\lambda + k\rho)}) P_\lambda(q^{2(\mu + k\rho)})$$

(here, as before, $[n] = \frac{q^n - q^{-n}}{q - q^{-1}}$).

Now, using the symmetry identity (4.3), let us replace $P_\lambda(q^{2(\mu + k\rho)})$ with $P_\mu(q^{2(\lambda + k\rho)}) \frac{g(\lambda)}{g(\mu)}$, where $g(\lambda) = \prod_{\alpha \in R^+} \prod_{i=0}^{k-1} [(\alpha, \mu + k\rho) + i]$. Then, after a short computation, we obtain the following recursive relation:

$$\sum_{\nu \in \Lambda_r : \mu + \nu \in P^+} \left( \prod_{\alpha \in R^+ : (\alpha, \nu) = -1} \frac{[(\alpha, \mu + k\rho) + k - 1][(\alpha, \mu + k\rho) - k]}{[(\alpha, \mu + k\rho)][(\alpha, \mu + k\rho) - 1]} \right) P_{\mu + \nu}(q^{2(\lambda + k\rho)})$$
$$(5.4) \qquad = X_r(q^{2(\lambda + k\rho)}) P_\mu(q^{2(\lambda + k\rho)}).$$

This relation is an equality between trigonometric polynomials in $\lambda$ satisfied for any dominant integral weight $\lambda$. Hence, it is satisfied identically, and we have

**Proposition 5.1.** *Macdonald's polynomials satisfy the recursive relations*
(5.5)
$$\sum_{\nu \in \Lambda_r : \mu+\nu \in P^+} \left( \prod_{\alpha \in R^+ : (\alpha,\nu)=-1} \frac{[(\alpha,\mu+k\rho)+k-1][(\alpha,\mu+k\rho)-k]}{[(\alpha,\mu+k\rho)][(\alpha,\mu+k\rho)-1]} \right) P_{\mu+\nu} = X_r P_\mu.$$

For example, if $n = 2, r = 1, \Lambda_r = \mathbb{C}_q^2$, relation (5.5) becomes the standard three-term recursive relation for the $q$-ultraspherical polynomials (formula (2.15) in [AI]).

*Remark.* Relations (5.5) determine the matrix of the operator of multiplication by $X_r$ in $\mathbb{C}[P]^W$ in the basis of Macdonald's polynomials.

**Corollary 5.2.** *The generalized characters $\chi_\lambda$ (see Section 2) satisfy*
(5.5)
$$\sum_{\nu \in \Lambda_r : \mu+\nu \in P^+} \left( \prod_{\alpha \in R^+ : (\alpha,\nu)=-1} \frac{[(\alpha,\mu+k\rho)+k-1][(\alpha,\mu+k\rho)-k]}{[(\alpha,\mu+k\rho)][(\alpha,\mu+k\rho)-1]} \right) \chi_{\mu+\nu} = X_r \chi_\mu.$$

Now let us assume that $\lambda$ is generic and consider the intertwining operator $\Phi_\lambda^k : M_\lambda \to M_\lambda \otimes U_{k-1}$ defined in Section 3. Let us introduce the $\psi$-functions

(5.7)
$$\psi_\lambda(x) = Tr|_{M_\lambda}(\Phi_\lambda^k x^h),$$

We would like to deduce recursive relations for the $\psi$-functions.

Looking at the expansions

(5.8)
$$\psi_{\lambda+(k-1)\rho}(x) = \sum_{\beta \in \lambda - Q^+} \psi_{\lambda+(k-1)\rho}^\beta x^{\lambda+(k-1)\rho-\beta},$$
$$\chi_\lambda(x) = \sum_{\beta \in \lambda - Q^+} \chi_\lambda^\beta x^{\lambda+(k-1)\rho-\beta},$$

we see that for each fixed $\beta$ the coefficients $\psi_{\lambda+(k-1)\rho}^\beta$, $\chi_\lambda^\beta$ coincide for sufficiently large $(\lambda, \rho)$, and that $\psi_{\lambda+(k-1)\rho}^\beta$ is a trigonometric rational function of $\lambda$. Therefore, Corollary 5.2 implies

**Proposition 5.3.** *The function $\psi_\lambda(x)$ satisfies the recurrence relations*
(5.5)
$$\sum_{\nu \in \Lambda_r : \mu+\nu \in P^+} \left( \prod_{\alpha \in R^+ : (\alpha,\nu)=-1} \frac{[(\alpha,\mu+\rho)+k-1][(\alpha,\mu+\rho)-k]}{[(\alpha,\mu+\rho)][(\alpha,\mu+\rho)-1]} \right) \psi_{\mu+\nu}(x) = X_r(x) \psi_\mu(x).$$

*This identity remains valid when the parameter $k$ is generic.*

## Appendix: ribbon graphs

For the sake of completeness we recall here the basic facts about the correspondence between ribbon graphs and representations of quantum groups, following as closely as possible [RT1, RT2].

A ribbon graph is an object in space formed by ribbons, which can be thought of as narrow strips of paper, and coupons, which are solid rectangles. Each ribbon and each coupon have a preferred direction. We require that each ribbon graph be located in between the horizontal planes $z = 0$ and $z = 1$, and the "free ends" of ribbons can be located either on the interval $[0, 1] \times 0 \times 0$ ("bottom") or on the interval $[0, 1] \times 0 \times 1$ ("top"). Some examples of ribbon graphs are shown on Figures 1 and 2. We always consider the ribbon graphs up to isotopy.

We consider "coloring" of ribbon graphs. That is, to each ribbon we assign a "color", i.e. an integral dominant weight $\lambda \in P^+$, and to each coupon we assign a homomorphism of $U_q\mathfrak{g}$-modules of the following type. Let us define the "bottom" and "top" sides of the coupon in such a way that the direction of the coupon is from top to bottom. If the colors of the ribbons with the ends on the "bottom" ("top") are $\lambda_1, \ldots, \lambda_k$ ($\mu_1, \ldots, \mu_m$, resp.), then coloring of the coupon is assigning to it a $U_q\mathfrak{g}$-intertwiner: $L_{\lambda_1}^{\varepsilon_1} \otimes \ldots \otimes L_{\lambda_k}^{\varepsilon_k} \to L_{\mu_1}^{\varepsilon_1} \otimes \ldots \otimes L_{\mu_m}^{\varepsilon_m}$, where $L^\varepsilon$ is either $L$ – if the directions agree – or $L^*$ – if they don't. For example, to the coupon shown on Figure 1b we must assign an intertwiner $L_\lambda \to L_\mu \otimes L_\nu^*$.

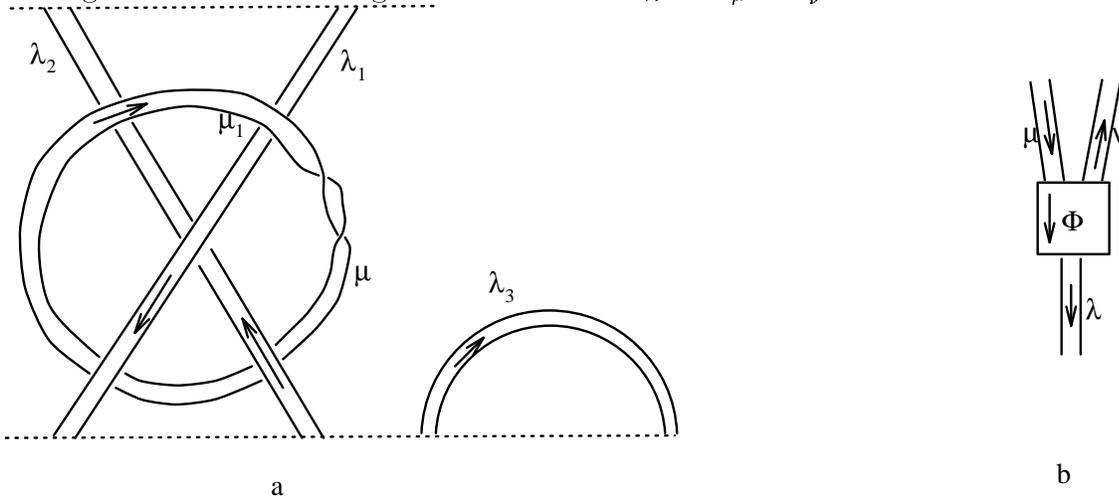

Figure 1

Such colored ribbon graphs can be multiplied in an obvious way if the directions and colors of the "bottom" of the first graph coincide with those of the "top" of the second one. We also have a notion of "tensor product" of ribbon graphs: if $\Gamma_1, \Gamma_2$ are ribbon graphs then $\Gamma_1 \otimes \Gamma_2$ is the ribbon graph which is obtained by placing $\Gamma_2$ to the right of $\Gamma_1$.

Then the main theorem, proved in [RT1] says that there is a unique way to assign to each colored ribbon graph $\Gamma$ a $U_q\mathfrak{g}$-homomorphism $F(\Gamma)$ so that the following conditions are satisfied:

a). if the colors of the ribbons with the ends on the "bottom" ("top") are $\lambda_1, \ldots, \lambda_k$ ($\mu_1, \ldots, \mu_m$, resp.), then $F(\Gamma)$ is an $U_q\mathfrak{g}$-intertwiner:

$$F(\Gamma): L_{\lambda_1}^{\varepsilon_1} \otimes \ldots \otimes L_{\lambda_k}^{\varepsilon_k} \to L_{\mu_1}^{\varepsilon_1} \otimes \ldots \otimes L_{\mu_m}^{\varepsilon_m},$$

where $L^\varepsilon$ is either $L$ or $L^*$ depending on the direction of the corresponding ribbon. For example, for the ribbon graph $\Gamma$ from the Figure 1, $F(\Gamma)$ is an intertwiner $L_{\lambda_1} \otimes L_{\lambda_2}^* \otimes L_{\lambda_3} \otimes L_{\lambda_3}^* \to L_{\lambda_2}^* \otimes L_{\lambda_1}$.

b). $F$ respects composition and tensor product.

c). the values of $F$ for the "elementary" graphs are shown below.

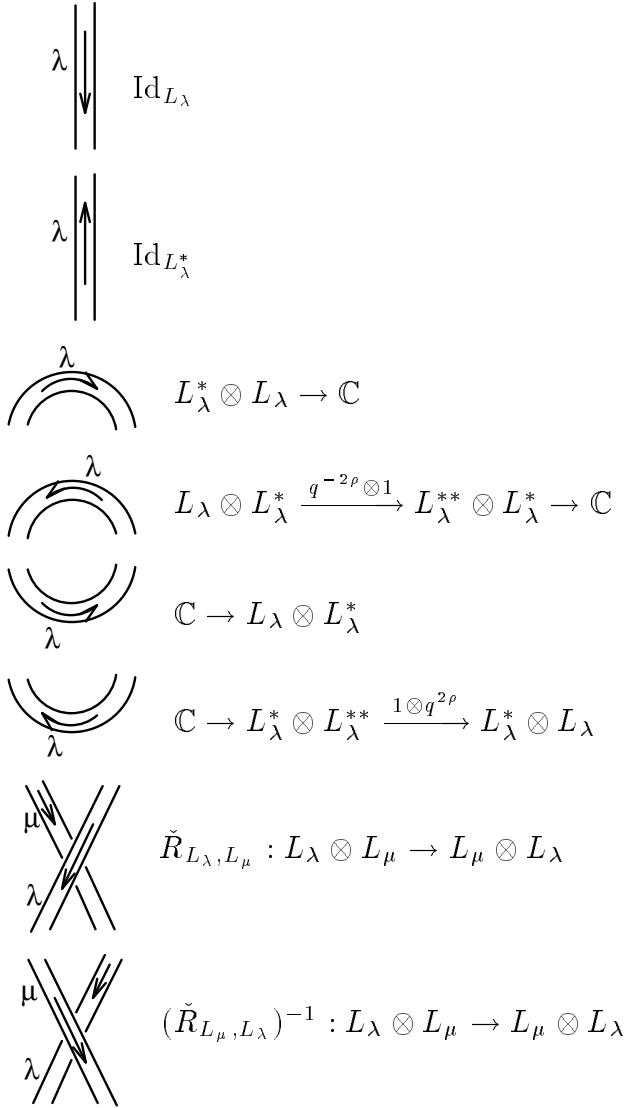

(*all unmarked arrows are canonical morphisms from* (1.1))

We will only need $F(\Gamma)$ for the graphs $\Gamma$ which have no twists. In this case, one can draw just the lines instead of the ribbons.

*Examples.*

1. If $f : L_{\lambda_1} \otimes \ldots \otimes L_{\lambda_k} \to L_{\lambda_1} \otimes \ldots \otimes L_{\lambda_k}$ is an intertwiner, then

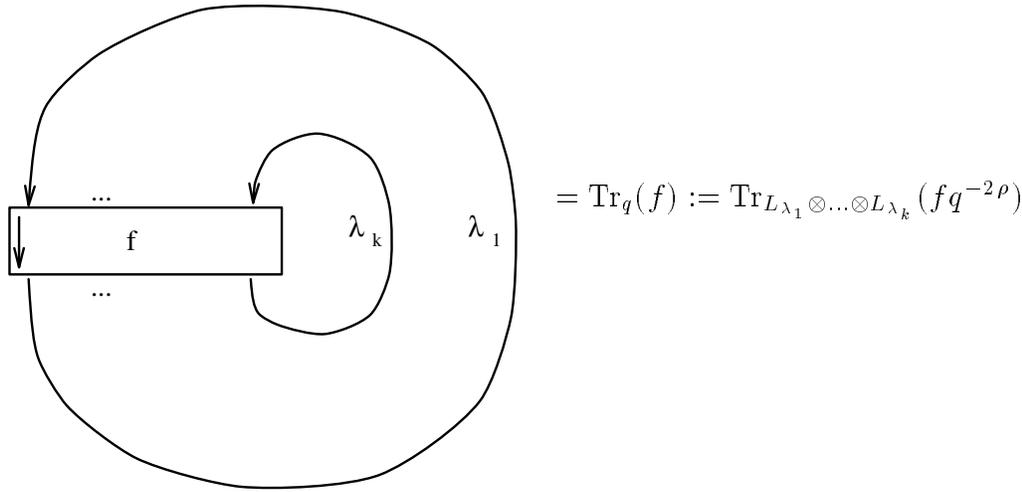

$$= \operatorname{Tr}_q(f) := \operatorname{Tr}_{L_{\lambda_1} \otimes \ldots \otimes L_{\lambda_k}}(f q^{-2\rho})$$

2. If $f : L_{\lambda_1} \otimes \ldots \otimes L_{\lambda_m} \to L_{\mu_1} \otimes \ldots \otimes L_{\mu_k}$ is an intertwiner, then

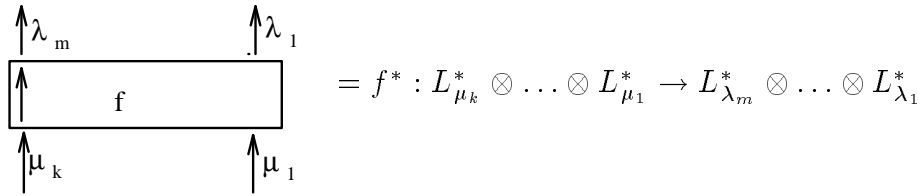

$$= f^* : L_{\mu_k}^* \otimes \ldots \otimes L_{\mu_1}^* \to L_{\lambda_m}^* \otimes \ldots \otimes L_{\lambda_1}^*$$

where $f^*$ is just the usual adjoint operator to $f$.